\documentclass{article}
\usepackage{beamerarticle}
\usepackage{beamerseminar}
\usepackage{beamertexpower}
\usepackage{hyperref}
\mode<presentation>
{
  \usetheme{}
  \setbeamercovered{transparent}
}
  
\usepackage[english]{babel}
\usepackage{xy}
\usepackage[T1]{fontenc} 
\usepackage[latin1]{inputenc}
\usepackage{amsfonts,amssymb,amsmath}
\usepackage{hyperref} 
\usepackage{graphics} 
\usepackage{multido, xcolor}
 
\newcommand{\BS}{\begin{frame}}
\newcommand{\ES}{\pause\end{frame}}

\newtheorem{theo}{\sf Theorem}

\newtheorem{coro}[theo]{\sf Corollary}
\newtheorem{defi}[theo]{\sf Definition}

\newtheorem{prop}[theorem]{\sf Proposition}

\newtheorem{exam}[theo]{\sc Example}

\def\qed{ \ \hfil$\square$}

\def\Ker{{\hskip0.3mm\rm  Ker\hskip0.5mm}}
\def\Im{{\hskip0.3mm\rm  Im\hskip0.5mm}}

\makeindex

\def \F {{\mathbb F}}

\def\ph{\varphi}

\def\End{{\rm End}}
\def\is{\buildrel \sim   \over \rightarrow}

\newcommand{\tr}{{\rm tr}}

\newcommand{\dd}{{\rm d}}

\newcommand{\Gal}{{\rm Gal}}

\newcommand{\GL}{{\rm GL}}
\newcommand{\GSp}{{\rm GSp}}

\newcommand{\Hom}{{\rm Hom}}

\newcommand{\ord}{{\rm ord}}
\newcommand{\SL}{{\rm SL}}

\newcommand{\A}{{\mathbb A}}
\newcommand{\C}{{\mathbb C}}
\newcommand{\N}{{\mathbb N}}
\newcommand{\Q}{{\mathbb Q}}
\newcommand{\R}{{\mathbb R}}
\newcommand{\Z}{{\mathbb Z}}
\newcommand{\Hb}{{\mathbb H}}

\newcommand{\Qb}{\overline{\mathbb Q}}

\newcommand{\chib}{\bar{\chi}}

\newcommand{\Ar}{{\mathcal A}}
\newcommand{\Ac}{{\mathcal A}}

\newcommand{\Cc}{{\mathcal C}}

\newcommand{\Dr}{{\mathcal D}}
\newcommand{\Hr}{{\mathcal H}}
\newcommand{\Mc}{{\mathcal M}}

\newcommand{\Lr}{{\mathcal L}}

\newcommand{\Oc}{{\mathcal O}}

\newcommand{\Sc}{{\mathcal S}}
\newcommand{\Tc}{{\mathcal T}}
\newcommand{\Pc}{{\mathcal P}}
\newcommand{\Qc}{{\mathcal Q}}

\newcommand{\Xc}{{\mathcal X}}
\newcommand{\Yc}{{\mathcal Y}}
\newcommand{\Zc}{{\mathcal Z}}

\newcommand{\al}{\alpha}
\newcommand{\De}{\Delta}
\newcommand{\de}{\delta}
\newcommand{\e}{\varepsilon}
\newcommand{\ep}{\varepsilon}
\newcommand{\Ga}{\Gamma}
\newcommand{\ga}{\gamma}

\newcommand{\La}{\Lambda}

\font\teneusm=eusm10 \font\seveneusm=eusm7 
\font\fiveeusm=eusm5 
\newfam\eusmfam 
\textfont\eusmfam=\teneusm 
\scriptfont\eusmfam=\seveneusm 
\scriptscriptfont\eusmfam=\fiveeusm

\def\nd{\not\hskip3.0pt\mid}
\def\nds{\not\hskip1.5pt\mid}
\def\mat #1,#2,#3,#4,{\left({#1\atop #3}{#2\atop #4}\right)}
\def\bra#1,{{\left\lbrace {#1}\right\rbrace}}
\def\ph{\varphi}

\def\si{\sigma}

\def \sgn {{\rm sgn}}

\def \diag {{\rm diag}}

\def \Spec {{\rm Spec}}

\def \Gal{{\rm Gal}}

\def \Id{{\rm Id}}
\def \Sp{{\rm Sp}}

\def\is{\buildrel \sim   \over \rightarrow}

\def\rr{\rangle}
\def\l1{\langle}

\newcommand{\B}{\left(\begin{array}{cc}}
\newcommand{\E}{\end{array}\right)}

\def \fns{{${}^{*)}$}}

\newcommand{\comm}[1]
{\fns\marginpar{$\boxed
{\hskip-6pt
{\small {\sf 
\begin{tabular} {l}
 #1
\end{tabular}
}
}
}
$
}
}
\def \?  {\comm{check ?}}

\usepackage{euscript}
\let\scr=\EuScript

\let\mathcal=\scr           
\def\ang#1,{{\left\langle {#1}\right\rangle}} 

\def\G{{\dbG}}

\def\qed{ \ \hfil$\square$}

\def\Ker{{\hskip0.3mm\rm  Ker\hskip0.5mm}}
\def\Im{{\hskip0.3mm\rm  Im\hskip0.5mm}}

\def \F {{\mathbb F}}

\def\ph{\varphi}

\def\End{{\rm End}}
\def\is{\buildrel \sim   \over \rightarrow}

\font\teneusm=eusm10 \font\seveneusm=eusm7 
\font\fiveeusm=eusm5 
\newfam\eusmfam 
\textfont\eusmfam=\teneusm 
\scriptfont\eusmfam=\seveneusm 
\scriptscriptfont\eusmfam=\fiveeusm

\font\tengothic=eufm10
\font\sevengothic=eufm7
\font\fivegothic=eufm5
\newfam\Gothic
\textfont\Gothic=\tengothic
\scriptfont\Gothic=\sevengothic
\scriptscriptfont\Gothic=\fivegothic
\def\go{\fam\Gothic\tengothic}
\def\a{{\mathfrak a}}
\def \c{{\go c}}
\def \m{{\go m}}
\def \e{{\go e}}

\def\nd{\not\hskip3.0pt\mid}
\def\nds{\not\hskip1.5pt\mid}
\def\mat #1,#2,#3,#4,{\left({#1\atop #3}{#2\atop #4}\right)}
\def\bra#1,{{\left\lbrace {#1}\right\rbrace}}
\def\ph{\varphi}

\def\si{\sigma}

\def\w{\omega}

\def \sgn {{\rm sgn}}

\def \diag {{\rm diag}}

\def \Re {{\rm Re}}

\def \Gal{{\rm Gal}}

\def \Id{{\rm Id}}
\def \Sp{{\rm Sp}}
\def\m{{\go m}}

\def\is{\buildrel \sim   \over \rightarrow}

\def\rr{\rangle}
\def\l1{\langle}
\def\Eq{\Longleftrightarrow}
\def\Numero{${\rm N\sp\circ}$}

\let\scr=\EuScript

\let\mathcal=\scr           
\let\db=\mathbb

\def\dbG{{\db G}}

\let\bold=\boldsymbol

\def\G{{\dbG}}

\def\vin{{ {\tiny \mid }  
\kern-7.29pt 
\bigcup }}

\def\ang#1,{{\left\langle {#1}\right\rangle}}

\def\m{{\frak m}}

\normalsize

\newcommand{\HH}{\mathbb H}

\def\ab{\A}

\def\qb{\Q}
\def\rb{\R}
\def\zb{\Z}

\def\dc{\Dr}
\def\hc{\Hr}

\def\mc{\Mr}

           \def\ab{{\Bbb A}}            
                      
           \def\cb{{\Bbb C}}            
\def\dc{{\cal D}}           

\def\hc{{\cal H}}

\def\mc{{\cal M}}

           \def\qb{{\Bbb Q}}            
           \def\rb{{\Bbb R}}

           \def\zb{{\Bbb Z}}

\def\diag {\mathop{\rm diag}\nolimits}
\def\SL {\mathop{\rm SL}\nolimits}
\def\Sp {\mathop{\rm Sp}\nolimits}

\font\tenrus=wncyr10
\font\sevenrus=wncyr7
\newfam\Rus
\textfont\Rus=\tenrus
\scriptfont\Rus=\sevenrus

\def\scr{\fam\Rus\sevenrus} 

\font\tenbrus=wncyb10
\font\eightbrus=wncyb8
\newfam\Brus
\textfont\Brus=\tenbrus
\scriptfont\Brus=\eightbrus

\font\tenirus=wncyi10
\font\eightirus=wncyi8
\newfam\Irus
\textfont\Irus=\tenirus
\scriptfont\Irus=\eightirus

\font\tenrus=wncyr10
\font\sevenrus=wncyr7
\newfam\Rus
\textfont\Rus=\tenrus
\scriptfont\Rus=\sevenrus

\def\scr{\fam\Rus\sevenrus} 
\input cyracc.def
\def\i1{\accent'044i}
\def\I1{\accent'044I}
\def\e1{\accent'040e}
\def\l1{l{}p1}

\def\qed{\quad\hbox{\hskip 1pt\vrule width 4pt height 6pt
          depth 1.5pt\hskip 1pt}}

\title[Graded structures and
quasimodular forms on classical groups
]
{Graded structures and differential operators on nearly holomorphic and
	quasimodular forms on classical groups}
\author{Alexei PANCHISHKIN
	\\ \small
	Institut Fourier, 
	Université Grenoble-1,  Université Grenoble Alpes\\
	100, rue des mathématiques, 38610, Gi\`eres, FRANCE}

\AtBeginSection[]
{
  \begin{frame}<beamer>
    \frametitle{} 
   \tableofcontents[currentsection]
\addtocounter{framenumber}{-1}
  \end{frame}
}

\ifarticle
\fi

\begin{document}
	
\date{}

\maketitle 
\begin{abstract}
	We wish to use graded structures \cite{KrVu87}, \cite{Vu01} on differential operators and
	quasimodular forms on classical groups and show that these structures provide a tool
	to construct $p$-adic measures and $p$-adic $L$-functions on the corresponding non-archimedean weight spaces.
	
	An approach to constructions of automorphic $L$-functions on unitary groups and their $p$-adic avatars is presented.
	For an algebraic group $G$ over a number field $K$ these $L$ functions 
	are  certain Euler products  $L(s,\pi, r, \chi)$.
	In particular, 
	our constructions cover the 
	$L$-functions in \cite{Shi00} via
	the doubling method of Piatetski-Shapiro and Rallis.		

	A $p$-adic analogue of $L(s,\pi, r, \chi)$ is a $p$-adic analytic function \\ $L_p(s,\pi, r, \chi)$ of $p$-adic arguments $s  \in \Z_p$, $\chi \bmod p^r$ which interpolates algebraic numbers  defined through the normalized critical values $L^*(s,\pi,r, \chi)$ of the corresponding complex analytic  $L$-function.
	We present a method using arithmetic nearly-holomorphic forms and general quasi-modular forms, related to algebraic automorphic forms.
	It gives a technique of constructing $p$-adic zeta-functions via 
	quasi-modular forms and their Fourier coefficients.
	
	Presented in a  talk for the
	INTERNATIONAL SCIENTIFIC CONFERENCE
	"GRADED STRUCTURES IN ALGEBRA AND THEIR APPLICATIONS"
	dedicated to the memory of Prof. Marc Krasner on
	Friday, September 23, 2016, 
	International University Centre (IUC), Dubrovnik, Croatia.

\end{abstract}


%
%
\setbeamertemplate{footline}
{%
	\begin{beamercolorbox}{section in foot} 
		\vskip-10pt
		\fbox{
			\insertpagenumber
		}
		\ {}
	\end{beamercolorbox}%
}

\section*{Introduction}
	Let $p$ be a prime number. Our purpose is to indicate a link of Krasner graded structures \cite{KrVu87}, \cite{Kr80}, \cite{Vu01} and constructions of  $p$-adic $L$-functions via distributions  and quasimodular forms on classical groups.

{ Krasner's graded structures are flexible and well adapted to various applications},  e.g. the rings and modules of differential  operators on classical groups and non-archimedean weight spaces.
	\subsection*{Arithmetical modular forms} 
		belong traditionally to the world of arithmetic, but also to the worlds of geometry, algebra and analysis.
		On the other hand,  it is very useful, to  attach  zeta-functions (or $L$-functions) to mathematical objects of various nature as certain generating functions (or as certain Euler products).
			\subsection*{Examples of modular forms,  quasimodular forms, 
zeta functions, and $L$-functions}
{\Large Eisenstein series\ }
	$\displaystyle E_k=1+\frac 2{\zeta(1-k)}\sum_{n= 1}^\infty\sum_{d|n}d^{k-1}q^n $ ($\in \Mc$, modular forms for even $k\ge 4, q=e^{2\pi i z}$), and $E_2\in \Qc\Mc $ {\it is  a quasimodular form}. 
	The ring of quasimodular forms, closed under differential operator ${\displaystyle D=q\frac{d}{dq}
		=\frac{1}{2\pi i}\frac{d}{dz}}$, used  in arithmetic, $\zeta(s)$ is the Riemann zeta function, $\displaystyle \zeta(-1)=-\frac 1{12}, \  E_2=1-24\sum_{n= 1}^\infty\sum_{d|n} d q^n$ which is also  a $p$-adic modular form (due to J.-P.Serre,\cite{Se73}, p.211)
		\\
		{\Large Elliptic curves} $E: y^2=x^3+ax+b$ , $a,b\in \Z$,  { A.Wiles's modular forms} $\displaystyle f_E=\sum_{n=1}^\infty a_nq^n$ with $a_p=p-Card E(\F_p)$ $(p\nd 4a^3+27b^2)$,
		 and the $L$-function $L(E,s)\displaystyle =\sum_{n=1}^\infty a_nn^{-s}$.

\subsection*{Zeta-functions or $L$-functions attached to various mathematical objects as certain Euler products.}

\begin{itemize}
\item $L$-functions link such objects to each other (a general form of functoriality); 

\item Special $L$-values answer fundamental questions about these objects in the form of a 
{ number (complex or $p$-adic)}.
\end{itemize}
Computing these numbers use integration theory
of Dirichlet-Hecke characters along  $p$-adic and complex valued measures.

This approach originates in the { Dirichlet class number formula} using the $L$-values in order  to compute class numbers of algebraic number fields through Dirichlet's $L$-series $L(s, \chi)$. 

The Millenium
{ BSD Conjecture} gives the rank of an elliptic curve $E$ as order of $L(E,s)$ at s=1
(i.e. the residue of its logarithmic derivative, see \cite{MaPa}, Ch.6).

		\
		
	
	
		\section{Graded  groups in the sense of Krasner}
\begin{defi} 
Let $G$ be a multiplicative group with the neutral element 1. A graduation of $G$ is a mapping $
\gamma: \Delta \rightarrow Sg(G), \delta \in \Delta
$ 
of a set $\Delta$ "the set of  grades of $\ga$ ") to the set  $Sg(G)$
 of subgroups of $G$ such that $G$ is the direct decomposition
 $$
 G=\bigoplus_{\delta\in \Delta}G_\delta, g=(g_{\delta})_{\delta\in \Delta}
 $$
The group $G$ endowed with such graduation is called {graded group}, 
and $g_{\delta}$ are { Krasner grade components} of $g$ with grade $\de$. 
\end{defi}

{\large Remarks} 
\begin{itemize}
\item The set $
H=\displaystyle\bigcup_{\delta\in \Delta}G_\delta
$
is called {\it the homogeneous part}  of $G$ for $\gamma$.
\item $x \in H$ are called  {\it homogeneous elements} of $G$ for $\gamma$
\item the elements $\delta \in \Delta$ are called {\it grades} of $\gamma$ 
\item  $\delta \in \Delta$ are called
{\it significant} (resp. empty) according as
$$
G_\delta\not =\{ 1\}, (resp. G_\delta =\{ 1\})
$$
and
$\Delta^*=\{ \delta\in \Delta ; G_\delta\not =\{ 1\} \}$ {\it the significative part} of $\Delta$. The graduation $\gamma$ is called {\it strict}  if $\Delta=\Delta^*$.

\end{itemize}
The choice of $\Delta$ is very flexible: we give examples of $p$-adic characters.	

\

\subsection{Graded  rings and modules}
	Moreover, M. Krasner introduced useful notions of graded  rings and modules.
	
	Let $(A; x+y, xy)$
	be a ring (not necessarily  associative) and let
	$\gamma:\Delta\to Sg(A; x+y)$
	be	a graduation of its addivive group. 
	The graduation $\gamma$ is called a graduation of a ring
	$(A; x+y, xy)= A$ if,
	for all 	$\xi, \eta\in \Delta$
	there exists a $\zeta\in \Delta$ such that	$A_\xi A_\eta\subset A_\zeta.$\\

	\subsection{Examples for $p$-adic groups $\Xc$,  and group rings} use the Tate field $\C_p=\hat {\bar \Q}_p$, the completion of ${\bar \Q}_p$, which is a fundamental object in  $p$-adic analysis, and thanks to Krasner we know that $\C_p$ is algebraically closed, see \cite{Am75}, Théorème
	2.7.1 ("Lemme de Krasner"), \cite{Kr74}. 
	This famous result allows to develop analytic functions and analytic spaces  over $\C_p$ (\cite{Kr74} , Tate, Berkovich...). and we  embed $incl_p : \bar\Q\hookrightarrow \bar\Q_p$.
		
	1) Algebraically, a $p$-adic measure $\mu$ on $\Xc$ is an element of the completed group ring	$A[[\Xc]]$, $A$ any $p$-adic subring of $\C_p$. 	\\	
	2)	The  $p$-adic $L$-function of $\mu$ is given by the evaluation  $\Lr_\mu(y):=y(\mu)$ on the
	group   $\Yc = \Hom_{cont}(\Xc,\C_p^*)$ of $\C_p^*$-valued characters  of $\Xc$.
	{\it The values $\Lr_\mu(y_j)$ on algebraic characters $ y_j\in \Yc^{alg}$ determine $\Lr_\mu$ iff they satisfy Kummer-type congruences}. 
	3) Our setting: a $p$-adic torus $T=\Xc$ of a unitary group $G$ attached to a CM field $K$ over $\Q$, a quadratic extension of a totally real field $F$, and an $n$-dimensional hermitian $K$-vector space $V$. Elements of  $\Yc^{alg}$ are identified with some algebraic characters  of the torus $T$ of the unitary group.
	
	\

\section{An extension problem.} From a subset $J=\Yc^{alg}$ of classical weights in  $\Yc=\Hom_{cont}(\Xc,\C_p^*)$, via the $A$-module $\Qc\Mc$ of quasimodular forms,  extend a mapping
	to the group ring $A[\Yc]$:
	\[
	\Lr: \Yc^{alg}\longrightarrow \mathop{
	\Qc\Mc}_{\displaystyle\circlearrowleft \Hr_A  }^{\displaystyle \curvearrowright
	\Dr_A } \longrightarrow \C_p, \ y_j \buildrel \ell \over \mapsto\Lr(y_j), y_j\in \Yc^{alg}
	\]
	($\Hr_A$ a Hecke\  algebra, $\Dr_A$ a ring of differential  operators over $A$),
	 or even to all continuous functions $\Cc(\Xc,\C_p)$, or just to local-analytic functions $\Cc^{loc-an}(\Xc,\C_p)$ such that the  
	values $\Lr(y_j)$ on $ y_j\in \Yc^{alg}$ are given algebraic $L$-values under the embedding
	$incl_p : \bar\Q\hookrightarrow \C_p$.

	Advantages of the $A$-module $\Qc\Mc$:
	
	1) nice Fourier expansion ($q$-expansions); 
	
	2) action of $D$ and of the ring of differential operators 
	$\Dr_A=A[D]$
	
	3) action of the Hecke algebra $\Hr_A$;
	
	4) projection $\pi_\alpha: \Qc\Mc^\alpha\to \Qc\Mc^\alpha$ to finite rank component  ("generalized eigenvectors of Atkin's $U$-operator") for any non-zero Hecke eigenvalue $\alpha$ of level $p$ ; $\ell$ goes through $\Qc\Mc^\alpha$ if $U^*(\ell)=\alpha^*\ell$.
	
{\it Solution} (extension of $\Lr$ to $\Cc(\Xc,\C_p)$) is given by the abstract Kummer-type congruences: 
$$
\forall x\in\Xc, \  \sum_j\beta_j y_j(x)\equiv 0 \mod p^N
\Longrightarrow
\sum_j\beta_j\ell (y_j)\equiv 0 \mod p^N  (\beta_j\in A).
$$
	These congruences imply the {\it  $p$-adic analytic continuation of the Riemann zeta function}:

	\

\subsection{Example: Mazur's $p$-adic integral
		and interpolation 
	}
	For any natural number $c > 1$ not divisible by $p$, there exists a $p$-adic measure 
	$\mu_c$ on $\Xc=\Z_p^*$, such that the special values
	$$
	\zeta(1-k)(1-p^{k-1})=
	{\int_{\Z_p^*}x^{k-1}d\mu_c(x) \over 1-c^{k}}
	\in \Q, (k \ge 2\mbox{ even })
	$$
	produce the {\it Kubota-Leopoldt $p$-adic zeta-function} 
	$
	\zeta_p: \Yc_p \to \cb_p
	$ 
	on the space $\Yc_p=Hom_{cont}(\Z_p^*, \cb_p^*)$ 
	as the  {\it $p$-adic Mellin transform} 
	$$
	\zeta_p(y)=
	{\int_{\Z_p^*}y(x)d\mu_c (x)\over 1-cy(c)}= {\Lr_{\mu_c}(y)\over 1-cy(c)}, 
	$$
	with a single simple pole at 
	$y=y_p^{-1}\in \Yc_p$, 
	where  
	$y_p(x)=x$ the inclusion character $\Z_p^* \hookrightarrow \C_p^*$ and $y(x)=\chi(x)x^{k-1}$  is a typical  arithmetical character ($y=y_p^{-1}$ becomes $k=0$, $s=1-k=1$). 
	
	{\large \it Explicitly}: Mazur's measure is given by $\mu_c(a+p^v\Z_p) = $\\
	$
	\frac {1}c\left[ \frac {ca}{p^v}\right] + 
	\frac {1-c}{2c}=\frac{1}c B_1(\{\frac {ca}{p^v}\})-B_1(\frac {a}{p^v})$,~$B_1(x)=x-\frac 1 2$,\\
	see \cite{LangMF}, Ch.XIII, 
$\ell$ is given by the {\it constant term of Eisenstein series}.	

	\
	
\subsection{Various exampes}
		The choice of $\Delta$ is quite flexible: we use the following examples: \\
		1) any direct decomposition
		$
		G=\bigoplus_{\delta\in \Delta}G_\delta
		$ gives a graduation where $\Delta$ is the set of homogeneous components and $\gamma=\Id : \Delta\to \Delta$,
		\\
		2) the regular representation of a commutative group algebra $A=K[G]$ with respect to characters of $G$, and of a subgroup 
$H\subset G$	gives a graduation of $A$\\
3) in particular, a Hodge structure of a complex vector space $V$ is a certain spectral decomposition for the real algebraic group ${\bf S}=\C^*$

Applications of the graduation 
$
G=\bigoplus_{\delta\in \Delta}G_\delta
$
are similar to the spectral theory, and the set
$\Delta$ extends the notion of the spectrum.
The homogeneous components $g_\delta$ in
$g=(g_\delta)_{\delta\in \Delta}$, are similar to the Fourier coefficients, they allow to extract essential information about $g\in G$.   
The same applies to {\it graded rings and modules} : $
A=\bigoplus_{\delta\in \Delta}A_\delta
$, $
M=\bigoplus_{\delta\in \Delta}M_\delta
$, the homogeneous Krasner grade components $a_\delta$,   $m_\delta$ simplify the $p$-adic constructions componentwise.   

	 \
	 

\subsection{Krasner grade components for proving Kummer-type congruences for $L$ and zeta-values}
		
	For more general  $L$-function $L(f,s)$,  of an automorphic form $f$ 
	one can prove certain Kummer-type congruences component-by-component using various Krasner grade components, with respect to {\it weights, Hecke-Dirichlet characters, eigenvalues of Hecke operators acting on spaces automorphic forms
	 (including Atkin-type $U_p$-operators), and the classical Fourier coefficients of quasimodular forms}.
	
	 It turns out that certain critical  $L$-values $L(f,s)$, expressed through Petersson-type product $\langle f, g_s \rangle $, reduces to $\langle \pi_\alpha(f), g_s \rangle $, where
	 $g_s$ is an explicit arithmetical automorphic form,  
	  $\al\ne 0$ is a  eigenvalue attached to $f$, and $\pi_\alpha(f)$ is the component given by the $\alpha$-characteristic projection, known to be in a fixed finite dimensional space (known for Siegel modular case, and extends to the unitary case).

		 \
		 

		
\section{Graded structures and $p$-adic measures.}
		For  a compact and totally disconnected topological	space  $X$, a $p$-adic ring $R$, let $C(X,R)$ be the $R$-module of continuous	maps  from $X$ to $R$ (with respect of the $p$-adic topology on $R$). For a  $M$ a $p$-adically complete $R$-module, an 
		{\it $M$-valued measure on $X$} is an element of the $R$-module
		$$
		Meas(X,M) = Hom_{\Z_p}(C(X,\Z_p),M)  \is Hom_R(C(X,R),M).
		$$

\subsection{An idea of construction of measures} uses certain   {\it  families  $\{ \phi_j\}_{j\in J}$ of test functions}, components of grade $j\in J$ in the group ring $R[Y] \subset C(X,R)$. In this case $J\subset \Delta=Y$.
A given family $m_j= \mu(\phi_j)\in M$ extends to  a measure $\mu$ on the $R$-linear span $R \langle\{ \phi_j\}_{j\in J}\rangle\subset C(X,R)$  
provided that certain Kummer-type congruences are satisfied  for $\mu(\phi_j)\in M$. 

Suppose $X$ is a profinite abelian group. Then $Meas(X,R)$ is identified with the completed
group ring $R[[X]]$ so that if a measure $\mu$ is identified with $f \in R[[X]]$, then for
any continuous character $\chi: X\to R^*_1$ with $R_1$ a $p$-adic $R$-algebra, $\mu(\chi) = \chi(f)$ (see	Section 6.1. Measures: generalities of \cite{EHLS}). 

\


\subsection{Example: $p$-adic characters of $\Z_p^*$} 
	This group $G=Hom_{cont}(\Z_p^*, \C_p^*)$ is used for constructions of $p$-adic zeta functions via	the $p$-adic integration theory
	 and the decomposition
	$Hom_{cont}(\Z_p^*, \C_p^*)=
	Hom_{cont}(\Gamma, \C_p^*)\times
	Hom_{cont}({\Z_p^*}^{tors}, \C_p^*)$, where 	
	$\Gamma=(1+p\Z_p)^*=\overline{\langle  1+p \rangle}$ for $p\ne 2$, when 
	$Hom_{cont}(\Gamma, \C_p^*)$  which is a $p$-adic analytic Lie group.

	{\it A graduation and an analytic structure on 
	$G$} are given by the above choice of the
	subgroup of tame characters $G_1= Hom_{cont}({\Z_p^*}^{tors}, \C_p^*)$ (a finite subgroup), and the
	subgroup of wild  characters $G_0=Hom_{cont}(\Gamma, \C_p^*)=U=\{1+t\in \C_p\ |t|_p<1\}$ (a $p$-adic disc).

	{\it The Mellin transform} of a $p$-adic distribution $\mu$ on $\Z^*_p$ gives an
	analytic function on the group of $p$-adic characters
	$$
	y \mapsto \Lr_\mu(y) =
	\int_{\Z^*_p}
	y(x)d\mu(y), \  y \in \Hom_{cont}(\Z^*_p, \C^*_p).
	$$
	According to {\it Iwasawa's theorem, there is one-to-one correspondence} between $p$-adic distributions $\mu$ on $\Z^*_p$ and bounded analytic functions on $\Yc$.
	
\

	\subsection{Methods of construction of $p$-adic $L$-functions}
	\begin{itemize}
			\item $p$-adic interpolation (starting from Kubota-Leopoldt in 1964);
		
		\item abstract Kummer congruences using a family of test elements (e.g. binomial coefficients as functions  with values in $\Z_p$, 
		or certain characters with values in $\bar\Q_p^*$, such as elements of $Y^{class}$ 
		
		\item the use of $p$-adic Cauchy integral (Shnirelman integral), e.g. a theorem of  Amice-Fresnel, see \cite{Ko80}, p.120. Let  
		$f(z) = \sum_n a_n z^n\in \C_p [[z]]$
		have the property that the coefficients $a_n$ can be $p$-adically
		interpolated i.e., there exists a continuous function $\phi: \Z_p\to \C_p$
		such that $\phi(n) = a_n$. Then $f$ (whose disc of convergence must be
		the open unit disc $D_0(1^-)$ is the restriction to $D_0 (1^- )$ of a Krasner analytic function $\tilde f$ on the complement of $D_1 (1^-)$, see also \cite{Kr74}.
	\item
	The Mellin transform of a $p$-adic measure $\mu$ on $\Xc$ gives an
	analytic function on the $p$-adic weight space $\Yc$, group of $p$-adic characters of $\Xc$:
	$$
	y \mapsto \Lr_\mu(y) =
	\int_{\Xc}
	y(x)d\mu(y), \  y \in \Yc=\Hom_{cont}(\Xc, \C^*_p).
	$$
	\
	
	\end{itemize}
	\section{Automorphic $L$-functions and their $p$-adic analogues} 
		{Our main objects  in this talk are automorphic $L$-functions and their $p$-adic analogues.}
		
		For an algebraic group $G$ over a number field $K$ these $L$ functions
		are defined as
		certain Euler products. 
		More precisely, we apply our constructions
		for the 
		$L$-functions studied in Shimura's book \cite{Shi00}.
		
		\begin{exam}
			[$G = \GL(2), K =\Q, L_f (s) =
			\sum_{n\ge 1}a_nn^{-s} , s \in \C$.]
			Here $f(z) =
			\sum_{n\ge 1}a_nq^n$ is a modular form on the upper-half plane
			$$
			H =\{ z \in \C, \Im(z) > 0\}= \SL(2)/SO(2),  q = e^{2\pi iz}.
			$$ 
			An Euler product has the form
			$$
			L_f (s) =
			\prod_{p \ primes}
			(1- a_pp^{-s} + \psi_f (p)p^{k-1-2s} )^{-1}
			$$
			where $k$ is the weight and $\psi_f$ the Dirichlet character of $f$. It is
			defined iff the automorphic representation $\pi_f$ attached to $f$ is irreducible, and
				$\pi_f$ is generated by the lift $\tilde f$ of $f$ to the group $G(\A)$. 
		\end{exam}
	

\

\subsection{A $p$-adic analogue of $L_f (s)$ (Manin-Mazur)}
It is a $p$-adic analytic function $L_{f ,p}(s, \chi)$ of $p$-adic arguments
$s \in\Z_p$,  $\chi \bmod p^r$ which interpolates algebraic numbers
$$
L^*_f (s,  \chi)/\omega^\pm\in\bar\Q\hookrightarrow\C_p = \hat{\bar\Q}_p \text{ (the Tate field)}
$$
for $1\le s\le k -1$, $\omega^\pm$ are periods of $f$ where the complex analytic
$L$ function of $f$ is defined for all $s\in\C$ so that in the absolutely
convergent case $\Re(s) > (k + 1)/2$,
$$
L_f^* (s,  \chi) = (2\pi)^{-s}\Gamma(s)
\sum_{n\ge 1}
\chi(n)a_nn^{-s}
$$
which extends to holomorphic function with a functional equation.
According to Manin and Shimura, this number is algebraic if the period $\omega^\pm$ is chosen according to the parity $\chi(-1)(-1)^{-s} = \pm 1$.

\subsection{Constructions of $p$-adic analogues
	of complex $L$-functions }
An irreducible automorphic representation $\pi$ of
adelic group $G(\A_K)$:
$$
L(s,  \pi,  r, \chi) =
\prod_{
	{\mathfrak p}_v\, primes\, in\, K}
\prod_{j=1}^m
(1-\beta _{j,{\mathfrak p}_v} N {\mathfrak p}_v^{-s} )^{-1}, 
\text{ where }\pi=\pi(f).
$$
 is an Euler product giving an $L$-function,
where $v \in \Sigma_K$ (places in
$K$), ${\mathfrak p} ={\mathfrak p}_v$,
$\alpha_{i,{\mathfrak p}}$ the Satake parameters of $\pi=\bigotimes_v\pi_v$, 
$$
\prod_{j=1}^m
(1-\beta _{j,{\mathfrak p}} X )
= \det(1_m-r(\diag(\alpha_{i,{\mathfrak p}})X)),
$$
$h_v =\diag(\alpha_{i,\mathfrak{p}})\in {}^LG(\C)$ (the Langlands group), $r : {}^LG(\C)\to \GL_m(\C)$ a finite dimensional representation, and
$\chi : \A^*_K/K^*\to \C^*$ is a character of finite order.

{\it Constructions extend to general automorphic representations on Shimura varieties} via the following tools:

\begin{itemize}
	\item Modular symbols and their higher analogues (linear forms on cohomology spaces related to automorphic forms)
	\item Petersson products with a fixed automorphic form, or
	\item linear forms coming from the Fourier coefficients (or Whittaker
	functions), or throught the
	\item CM-values (special points on Shimura varieties),
\end{itemize}

\

\subsection{Accessible cases: symplectic and unitary groups}
	\begin{itemize}
		\item
		$G = \GL_1$ over $\Q$ (Kubota-Leopoldt-Mazur) for the Dirichlet
		$L$-function $L(s, \chi)$.
		\item
		$G = \GL_1$ over a totally real field $F$ (Deligne-Ribet, using algebraicity result by Klingen).
		
		\item
		$G = \GL_1$ over a CM-field $K$, i.e. a totally imaginary extension of a totally real field $F$ (N.Katz, Manin-Vishik).
		\item {the Siegel modular case}
		$G = \GSp_n$ (the Siegel modular case, $F=\Q$). 
		\item{General symplectic and unitary groups  over a CM-field $K$}
		Certain Euler products in Chapter 5 of \cite{Shi00}, 
		with critical values computed in Chapter 7, Theorem 28.8 using general nearly holomorphic arithmetical automorphic  forms for the group
		$$
		G= G(\varphi) =\{ \al\in \GL_{n+m}(K)| \al\varphi\al^*=\nu(\alpha)\varphi
		\},\nu(\al)\in F^*,
		$$
		where
		$
		\varphi=
		\begin{pmatrix}
		1_n&0\\
		0&1_m
		\end{pmatrix}
		$, or
		$\varphi=\eta_n =
		\begin{pmatrix}
		0&1_n\\
		-1_n& 0
		\end{pmatrix}, n=m 
		$, see also
		Ch.Skinner and E.Urban \cite{MC} and Shimura G., \cite{Shi00}.
	\end{itemize}
	
\

\subsection{Nearly holomorphic modular forms for $\GL(2)$} in N.M.Katz \cite{Ka76} 
		 used real-analytic and quasimodular forms coming
	from derivatives of holomorphic forms, and $p$-adic modular forms. 
	
	A relation real-analytic $\leftrightarrow p$-adic modular forms comes from the
	notion of $p$-adic modular forms invented by J.-P.Serre \cite{Se73} as
	$p$-adic limits of $q$-expansions of modular forms with rational
	coefficients for $\Gamma= \SL_2(\Z)$, using the ring $M_p$ of $p$-adic modular forms contains $\Mc = \oplus_{k\ge 0}\Mc_k(\Gamma, \Z) = \Z[E_4,  E_6]$, and it contains $E_2$ as an element with the $q$-expansion
	$E_2 = 1- 24
	\sum_{n\ge 1}\sigma_1(n)q^n$. 
	On the other hand, the function of $z=x+iy$,  
	$$
	\tilde E_2 = -\frac {3}{\pi y}
	+ E_2 =-12S + E_2,  \text{ where }
	S =
	\frac {1}{4\pi y},
	$$
	is a nearly holomorphic modular form (that is,  its coefficients are
	polynomials of $S$ over $\Q$). 
Then ${\mathcal N}=\Z[\tilde E_2, E_4, E_6]$ is the ring of such forms, 	
	$\tilde E_2|_{S=0} = E_2$
	and $E_2$ is a $p$-adic modular form. 
	
	Elements of the ring $\Qc{\Mc}={\mathcal N}|_{S=0}
	=\Z[E_2, E_4, E_6]$ are called quasimodular forms hence such forms are all $p$-adic modular forms. 
	These
	phenomena are quite general and used in computations and
	proofs.
	S.Boecherer extended these
	results to the Siegel modular case.	
	
	\
	

	\section{Automorphic $L$-functions attached to symplectic and unitary groups}
	Let us briefly describe the $L$-functions attached to symplectic and unitary groups as certain Euler products in Chapter 5 of \cite{Shi00}, 
	with critical values computed in Chapter 7, Theorem 28.8 using general nearly holomorphic arithmetical automorphic  forms for the group
	$$
	G= G(\varphi) =\{ \al\in \GL_{n+m}(K)| \al\varphi\al^*=\nu(\alpha)\varphi
	\},\nu(\al)\in F^*,
	$$
	where
	$
	\varphi=
	\begin{pmatrix}
	1_n&0\\
	0&1_m
	\end{pmatrix}
	$, or
	$\varphi=\eta_n =
	\begin{pmatrix}
	0&1_n\\
	-1_n& 0
	\end{pmatrix}, n=m 
	$, see also
	Ch.Skinner and E.Urban \cite{MC} and Shimura G., \cite{Shi00}.
\subsection{The groups and automorphic forms in \cite{Shi00}}
Let $F$ be a totally real algebraic number field, $K$ be a totally
imaginary quadratic extension of $F$ and $\rho$ be the generator of
$\Gal(K/F)$. 
Take $\eta_n =
\begin{pmatrix}
0&1_n\\
-1_n& 0
\end{pmatrix}
$
and define
\begin{align*} &G = \Sp(n, F) &&\text{(Case Sp)}
\\ & G = \{\al\in \GL_{2n}(K)| \al\eta_n\al^*= \eta_n\} &&\text{(Case UT = unitary tube)}
\\ &G = \{\al\in \GL_{2n}(K)| \al T\al^* = T\} &&\text{(Case UB =
	unitary ball)}
\end{align*}
according to three cases.
Assume $F=\Q$ for a while. The group of the real points $G_\infty$ acts on the associated domain
\begin{align*}
{\mathcal H} =
\begin{cases}\{ z \in M(n, n,\C)\  | \ {}^tz = z, \Im(z) > 0\} & \text{(Case Sp)}\\
\{ z \in M(n, n,\C)\  | \  i(z^*-z) > 0\} & \text{(Case UT)}\\
\{ z \in M(p, q,\C)\  | \  1_q-z^*z > 0\} & \text{(Case UB).}
\end{cases}
\end{align*}
$(p, q), p + q = n$, being the signature of $iT$ . Here $z^* = {}^t\bar z$ and >
means that a hermitian matrix is positive definite.
In Case UB, there is the standard automorphic factor
$M(g,  z)$, $g\in G_\infty$,  $z \in \Hr$ taking values in 
$\GL_p(\C) \times\GL_q(\C)$.

\

	\subsection{Shimura's arithmeticity  \cite{Shi00},  $p$-adic zeta functions
		and nearly-holomorphic forms on classical groups}
	Automorphic $L$-functions and their $p$-adic analogues can be obtained
	for quite general automorphic representations on Shimura varieties by constructing $p$-adic distributions out of algebraic numbers attached to automorphic forms. 
	This means that these numbers satisfy certain Kummer-type congruences established 
	in different  ways: via
	\begin{itemize}
		\item Normalized Petersson products with a fixed automorphic form, or
		\item linear forms coming from the Fourier coefficients (or Whittaker
		functions), or through the
		\item CM-values (special points on Shimura varieties),
		see The Iwasawa Main Conjecture for GL(2) by C. Skinner and E.
		Urban, [MC], Shimura G., Arithmeticity in the theory of automorphic forms \cite{Shi00}.

	\end{itemize}
	
	The combinatorial structure of the Fourier coefficients of the
	holomorphic forms used in these constructions is quite complicated.
	
	In order to prove the congruences needed for the $p$-adic constructions, we use a simplification due to nearly-holomorphic and general
	quasi-modular forms, related to algebraic automorphic forms.
	In this paper, a new method of constructing $p$-adic zeta-functios is
	presented using general quasi-modular forms and their Fourier	coefficients.

\	

\subsection{Algebraicity and congruences of the critical values of the zeta functions} of automorphic forms on unitary and
	symplectic groups, we follow the review by H.Yoshida \cite{YS} of Shimura's book \cite{Shi00}.
	
	Shimura's mathematics developed by stages:
	
	(A) Complex multiplication of abelian varieties =>
	
	(B) The theory of canonical models = Shimura varieties =>
	
	(C) Critical values of zeta functions and periods of automorphic forms.
	
	(B) includes (A) as 0-dimensional special case of canonical models.
	The relation of (B) and (C) is more involved, but (B) provides a solid foundation of the notion of the arithmetic automorphic forms.
	Also unitary Shimura varieties have recently attracted much interest
	(in particular by Ch. Skinner and E. Urban), see \cite{MC}, in relation
	with the proof of the The Iwasawa Main Conjecture for $\GL(2)$.
	
	\
	
\subsection{Integral representations and critical values of the zeta functions}
	In	Cases Sp and UT, Eisenstein series $E(z,  s)$ associated to the	maximal parabolic subgroup of $G$ of Siegel type is introduced. 
		Its analytic behavior and those values of $\si\in 2^{-1}\Z$ at which $E(z,  \sigma)$ is nearly
		holomorphic and arithmetic are studied in \cite{Shi00}. 
		This is achieved by proving
		a relation giving passage from $s$ to $s - 1$ for $E(z, s)$, involving a
		differential operator, then examining Fourier coefficients of
		Eisenstein series using the theory of confluent hypergeometric functions on
		tube domains.
		
		For a Hecke eigenform $f$ on $G_\A$ and an algebraic Hecke character $\chi$
		on the idele group of $K$ (in Case Sp, $K = F$), the zeta function
		{\it ${\mathcal Z}(s,  f,  \chi)$,
			 an Euler product extended over prime ideals of $F$, the degree of the Euler factor is $2n + 1$ in Case
		Sp, $4n$ in Case UT, and $2n$ in Case UB}, except for finitely many prime ideals, see  Theorem 19.8 on Euler products of Chapter 5, \cite{Shi00}.
		
		This zeta function is almost the same as the so called standard $L$-function attached to $f$ twisted by $\chi$ but it turns out to be more general in the unitary case, see also \cite{EHLS}.
		
		Main results on critical values of 
		the $L$-functions studied in Shimura's book \cite{Shi00}
		is stated in Theorem 28.5, 28.8 (Cases Sp, UT), and in Theorem 29.5 in Case UB.
		
		\
		
\subsection{Shimura's Theorem: algebraicity of critical values in Cases Sp and UT, p.234 of \cite{Shi00}}
			Let ${\bf f}\in {\mathcal V} (\bar\Q)$ be a non zero  arithmetical automorphic form of type Sp or UT. 
			Let $\chi$ be a Hecke character of $K$ such that 
			$\chi_{\bf a}(x)=x_{\bf a}^\ell |x_{\bf a}|^{-\ell}$ with $\ell\in\Z^{{\bf a}}$, \text{ and let }$\sigma_0\in 2^{-1}\Z$.
			Assume, in the notations of Chapter 7 of \cite{Shi00} onr the weights $k_v, \mu_v, \ell_v$, that 
			\begin{align*}
			&\text{\rm Case Sp} &&
			\  2n+1-k_v+\mu_v \le 2\si_0\le k_v-\mu_v,  
			\\&{}  && \  \text{ where }
			\mu_v=0  \text{ if } [k_v]-l_v \in 2\Z \\&{}  && \  \text{ and } \mu_v=1
			 \text{ if } [k_v]-l_v \not\in 2\Z; \ 
			\si_0-k_v+\mu_v
			\\&{}  && \  \text{ for every } v \in  {\bf a}  \text{ if }\si_0>n
			\text{ and } 
			\\&{}  &&\ \si_0-1-k_v+\mu_v \in 2\Z\text{ for every } v \in  {\bf a}
			\text{ if }\si_0\le n. \\
			&\text{\rm Case UT} &&\  4n-(2k_{v\rho}+\ell_v)\le 2\si_0\le 
			m_v-| k_v-k_{v\rho}-\ell_v |\\&{}  && \  \text{ and } 2\si_0-\ell_v\in 2\Z  \text{ for every  } v \in  {\bf a}.
			\end{align*}
			Further exclude the following cases
			\begin{align*}
			&\text{\rm (A) Case Sp} && \si_0=n+1,F=\Q  
			\text{ and }
			\chi^2=1;
			\\ &\text{\rm (B) Case Sp} &&\si_0=n+(3/2),F=\Q; \chi^2=1  
			\text{ and }
			[k]-\ell\in 2\Z
			\\ &\text{\rm (C) Case Sp} &&\si_0=0,
			{\mathfrak c}={\mathfrak g} 
			\text{ and }
			\chi=1;
			\\ &\text{\rm (D) Case Sp} &&0<\si_0\le n, {\mathfrak c}={\mathfrak g}, \chi^2=1
			\text{ and the conductor of } \chi \text{ is }{\mathfrak g};
			\\&\text{\rm (E) Case UT} && 2\si_0=2n+1,F=\Q, \chi_1=\theta,  
			\text{ and }
			k_v-k_{v\rho}=\ell_v;
			\\&\text{\rm (F) Case UT} && 0\le 2\si_0 < 2n, {\mathfrak c}={\mathfrak g},
			\chi_1=\theta^{2\si_0} 	\text{ and the conductor of } \chi 	\text{ is } {\mathfrak r}
			\end{align*}
			Then 
			$$
			{\mathcal Z}(\si_0, {\bf f}, \chi)/\langle  {\bf f},  {\bf f} \rangle\in\pi^{n|m|+d\varepsilon}\bar \Q, 
			$$
			where $d = [F:\Q]$, $|m|=\sum_{v\in {\bf a}}m_v$, and
			$$\varepsilon=
			\begin{cases}
			(n+1)\si_0-n^2-n,&\text{Case Sp}, k\in \Z^{{\bf a}}, 
			\text{ and } \si_0>n_0),
			\\
			n\si_0-n^2, &\text{Case Sp}, k \not\in \Z^{{\bf a}}, 
			\text{or} \si_0\le n_0),
			\\ 2n\si_0-2n^2+n
			&\text{Case UT}
			\end{cases}
			$$
Notice that $\pi^{n|m|+d\varepsilon}\in\Z$ in all cases; if $k\not \in\Z^{\bf a}$,  the above parity condition on $\sigma_0$ shows that  $\sigma_0+k_v\in\Z$, so that $n|m|+d\varepsilon\in\Z$.

\
		
	\subsection{A $p$-adic analogue of Shimura's Theorem (Cases Sp and UT)}
		represent algebraic parts of critical values as values of certain $p$-adic analytic zeta function, Mellin transform  $y\mapsto \Lr_\mu(y)$ of a $p$-adic distributions $\mu$ on $\Xc=T$, an
		analytic function on the group $\Yc$ as above.
		
		\
		The construction of  $\mu$ uses congruences of Kummer type between the Fourier coefficients of modular forms.	Suppose that we are given some $L$-function ${\Zc}^*(s, {\bf f}, \chi)$ attached to a quasimodular form ${\bf f}$ and assume that for infinitely many "critical pairs" $(s_j ; \chi_j )$ one has an integral representation
		$
		{\Zc}^*(s_j, {\bf f}, \chi_j)
		=\langle {\bf f} , h_j \rangle
		$
		with all $h_j =
		\sum_\Tc b_{j, \Tc} q^\Tc \in \Mc$
		in a certain
		finite-dimensional space $\Mc$ containing ${\bf f}$ and defined over $\bar{\Q}$.
		We want to prove the following Kummer-type congruences:
		$$
		\forall x\in\Z^*_p, \  \sum_j\beta_j\chi_jx^{k_j}\equiv 0 \mod p^N
		\Longrightarrow
		\sum_j\beta_j\frac{ {\Zc}^*(s, {\bf f}, \chi)}{\langle {\bf f} , {\bf f} \rangle}\equiv 0 \mod p^N
		$$
		$$
		\beta_j \in \bar \Q, k_j =
		\begin{cases}
		s_j -s_0 & \text{ if }s_0 = \min_j s_j \text{ or }\\
		s_0-s_j&\text{ if }s_0 =  \max_j s_j .
		\end{cases}
		$$
		{Computing the Petersson products}
		of a given quasimodular form
		${\bf f} (Z) =\sum_\Tc a_\Tc q^\Tc\in \Mc_\rho(\bar{\Q})$
		by another quasimodular form
		$h (Z) =\sum_\Tc b_\Tc q^\Tc\in \Mc_{\rho^*}(\bar{\Q})$ uses a linear form
		$$
		\ell_{\bf f} : h \mapsto \frac{\langle {\bf f} , h \rangle}{\langle {\bf f} , {\bf f} \rangle}
		\text{ defined over a subfield }k \subset \bar\Q.
$$

	\subsection{Using algebraic and $p$-adic modular forms}
		There are several methods to compute various $L$-values starting
		from the constant term of the Eisenstein series in \cite{Se73},
		$$
		G_k (z) =
		\frac{\zeta(1-k)} 2
		+
		\sum_{
			n=1}^\infty
		\si _{k-1}(n)q^n =
		\frac{\Gamma(k)}
		{(2\pi i )^k}
		\sum_{(c,d)} {}'
		(cz + d)^{-k} \ \ (\text{for} k\ge 4),
		$$
		and using Petersson products of nearly-holomorphic Siegel modular
		forms and arithmetical automorphic forms as in \cite{Shi00}:
		
		the Rankin-Selberg method,
		
		the doubling method (pull-back method).
		
		A known example is the standard zeta function 
		$D(s,  f, \chi)$ of a
		Siegel cusp eigenform   $f\in \Sc_n^k (\Gamma)$ of genus $n$ (with local factors of
		degree 2$n$ + 1) and $\chi$ a Dirichlet character.
		
		{\bf Theorem} (the case of even genus $n$ (\cite{Pa91}, \cite{CourPa}), via the
		Rankin-Selberg method) gives a $p$-adic interpolation of the
		normailzed critical values $D^*(s, f, \chi)$ using Andrianov-Kalinin integral representation of these values $1 + n -k \le  s \le  k - n$
		through the Petersson product $\langle f,  \theta_{T_0}\delta^rE\rangle$ where $\delta^r$ is a certain
		composition of Maass-Shimura differential operators, $\theta_{T_0}$ a theta-series of weight $n/2$, attached to a fixed $n \times n$ matrix $T_0$.

\

	\subsection{Using $p$-adic doubling method}		
		{\bf Theorem}
			($p$-adic interpolation of
			$D(s, f, \chi)$)\\
			 (1) (the case of odd genus (Boecherer-Schmidt, \cite{BS00}) 
			\ \\
		{\em	Assume that $n$ is arbitrary genus, and a prime $p$ ordinary then there exists a $p$-adic interpolation of
			$D(s, f, \chi)$ 
			
			{\rm (2)} (Anh-Tuan Do (non-ordinary case, PhD Thesis, 2014)), via
				the doubling method)}
			\ \\  Assume that $n$ is arbitrary genus, and  $p$ an arbitrary prime not dividing level of $f$ 
			then there exists a $p$-adic interpolation of
			$D(s, f, \chi)$. 	
			
		
		Proof
		uses the following Boecherer-Garrett-Shimura identity
		(a pull-back formula, related to the Basic Identity of {P}iatetski-Shapiro and Rallis, \cite{GPSR}),
		which allows to compute the critical values through certain double
		Petersson product by integrating over 
		$z\in \HH_n$ the identity:
		$$
		\Lambda(l + 2s, \chi)D(l + 2s-n,  f, \chi) f =
		\langle f (w),  E^{2n}_{l, \nu, \chi, s} (\diag[z, w])_w, 
		$$
		Here $k = l +\nu, \nu \ge 0, \Lambda(l + 2s, \chi)$ is a product of special values of
		Dirichlet $L$-functions and $\Gamma$-functions,
		$E^{2n}_{l, \nu, \chi, s}$ a higher twist of a Siegel-Eisenstein series on $(z,w)\in \Hb_n\times \Hb_n$ (see \cite{Boe85},
		\cite{BS00}).

\

\subsection{Injecting nearly-holomorphic forms into $p$-adic modular forms}
	A recent discovery by Takashi Ichikawa (Saga University),
	[Ich12], J. reine angew. Math., [Ich13]
	allows to inject nearly-holomorphic arithmetical (vector valued)
	Siegel modular forms into $p$-adic modular forms.
	Via the Fourier expansions, the image of this injection is
	represented by certain quasi-modular holomorphic forms like
	$E_2 = 1- 24
	\sum_{n\ge 1}
	\si_1(n)q^n$, with algebraic Fourier expansions.
	This description provides many advantages, both computational
	and theoretical, in the study of algebraic parts of Petersson
	products and $L$-values, which we would like to develop here.
	In fact,
	the realization of nearly holomorphic forms as $p$-adic modular forms has been
	studied by Eric Urban, who calls them "Nearly overconvergent modular forms" \cite{U14}, Chapter 10.
	
	Urban only treats the elliptic modular case in that paper, but I believe
	he and Skinner are working on applications of a more general theory.
	This work is related to a recent preprint \cite{BoeNa13} by S. Boecherer
	and Shoyu Nagaoka where it is shown that Siegel modular forms of
	level $\Gamma_0(p^m)$ are $p$-adic modular forms. Moreover they show that
	derivatives of such Siegel modular forms are $p$-adic. Parts of these
	results are also valid for vector-valued modular forms.
	\

	\subsection{Arithmetical nearly-holomorphic Siegel modular forms}
	Nearly-holomorphic Siegel modular forms over a subfield $k$ of $\C$ are
	certain $\C^d$ -valued smooth functions $f$ of 
	$Z = X + \sqrt{-1}Y  \in \Hb_n$
	given by the following expression
	$
	f (Z) =
	\sum_T
	P_T (S)q^T 
	$
	where $T$ runs through the set $B_n$ of all half-integral semi-positive
	matricies, $S = (4\pi Y )^{-1}$ a symmetric matrix,
	$q^T = \exp(2\pi \sqrt{-1} \tr(TZ)), P_T (S)$ are vectors of degree d whose
	entries are polynomials over $k$ of the entries of $S$.
	
	{\it A geometric construction of such arithmetical forms uses the algebraic theory af moduli spaces of abelian varieties.}
	Following \cite{Ha81}, consider the columns $Z_1, Z_2, \dots Z_n$ of $Z\in \Hb_n$ and the
	$\Z$-lattice $L_Z$ in $\C^n$ generated by $\{
	E_1, E_2, \dots E_n, Z_1, Z_2, \dots Z_n \}, 
	$
	where
	$E_1, E_2, \dots E_n$ are the columns of the identity matrix $E$. 
	The torus $A_Z = \C^n/L_Z$ is an abelian variety, and there is an analytic family
	$\Ar \to \Hb_n$ whose fiber over the point $Z\in \Hb_n$ is $A_Z$ .
	Let us consider the quotient space $\Hb_n/\Gamma(N)$ of the Siegel upper half space $\Hb_n$ of degree $n$ by the integral symplectic group
	$$
	\Gamma(N) =
	\left\{
	\gamma=
	\begin{pmatrix}
	&A_\gamma 
	&B_\gamma
	\\
	&C_\gamma 
	&D_\gamma
	\end{pmatrix} 
	\ \Big| \
	\begin{matrix}
	& A_\gamma \equiv D_\gamma \equiv 1_n \\
	&B_\gamma \equiv C_\gamma \equiv 0_n
	\end{matrix}
	\right\}, 
	$$
	giving the moduli space classifying
	principally polarized abelian schemes of relative dimension $n$ with a symplectic level $N$ structure.

\
	
	\subsection{Applications to constructions of $p$-adic $L$-functions}
	We present here a survey
	of some methods of construction of $p$-adic $L$-functions. 
	Two important ideas
	that are not as well known as they should be are developed biefly in this section.
	
	There exist two kinds of $L$-functions
	\begin{itemize}
		\item Complex $L$-functions (Euler products) on 
		$\C = \Hom(\R^*_+; \C^*)$.
		\item $p$-adic $L$-functions on the $\C_p$-analytic group 
		$\Hom_{cont}(\Z^*_p, \C_p^*)$
		(Mellin transforms $\Lr_{\mu}$ of $p$-adic measures $\mu$ on $\Z^*_p$).
	\end{itemize}
	
	Both are used in order to obtain a number ($L$-value) from an
	automorphic form. Such a number can be algebraic (after
	normalization) via the embeddings,
	$$
	\bar\Q \hookrightarrow \C, \bar\Q \hookrightarrow \C_p=\hat{\bar\Q}_p 
	$$  
	and we may compare the complex and $p$-adic $L$-values at many
	points.
	
	\
	
\subsection*{ How to define and to compute $p$-adic $L$-functions?}
	The Mellin transform of a $p$-adic distribution $\mu$ on $\Z^*_p$ gives an
	analytic function on the group of $p$-adic characters
	$$
	y \mapsto \Lr_\mu(y) =
	\int_{\Z^*_p}
	y(x)d\mu(x), \  y \in \Hom_{cont}(\Z^*_p, \C^*_p)
	$$
	A general idea is to construct $p$-adic measures directly from Fourier
	coefficients of modular forms proving Kummer-type congruences for
	$L$-values. 
	Here we present a new method to construct $p$-adic
	$L$-functions via quasimodular forms.
	
	\
	
\subsection*{How to prove Kummer-type congruences
		 using the Fourier coefficients?}

	Suppose that we are given some $L$-function $L^*_f (s, \chi)$ attached to a Siegel modular form $f$ and assume that for infinitely many "critical pairs" $(s_j ; \chi_j )$ one has an integral representation
	$
	L^*_f (s,  \chi)=\langle f , h_j \rangle
	$
	with all $h_j =
	\sum_\Tc b_{j, \Tc} q^\Tc \in \Mc$
	in a certain
	finite-dimensional space $\Mc$ containing $f$ and defined over $\bar{\Q}$.
	We want to prove the following Kummer-type congruences:
	$$
	\forall x\in\Z^*_p, \  \sum_j\beta_j\chi_jx^{k_j}\equiv 0 \mod p^N
	\Longrightarrow
	\sum_j\beta_j\frac{ L^*_f (s_j,  \chi)_j}{\langle f , f \rangle}\equiv 0 \mod p^N
	$$
	$$
	\beta_j \in \bar \Q, k_j =
	\begin{cases}
	s_j -s_0 & \text{ if }s_0 = \min_j s_j \text{ or }\\
	s_0-s_j&\text{ if }s_0 =  \max_j s_j .
	\end{cases}
	$$
{Computing the Petersson products}
of a given modular form
$f (Z) =\sum_\Tc a_\Tc q^\Tc\in \Mc_\rho(\bar{\Q})$
by another modular form
$h (Z) =\sum_\Tc b_\Tc q^\Tc\in \Mc_{\rho^*}(\bar{\Q})$ uses a linear form
$$
\ell_f : h \mapsto \frac{\langle f , h \rangle}{\langle f , f \rangle}
$$
defined over a subfield $k \subset \bar\Q$.

	\
	
	\subsection{Proof of the main congruences}
Thus the Petersson product in $\ell_{\bf f}$ can be expressed through
the Fourier coeffcients of $h$ in the case when there is a finite basis
of the dual space consisting of certain Fourier coeffcients:
$\ell_{\Tc_i} : h \mapsto b_{\Tc_i} (i = 1, \dots, n)$.
It follows that 
$\ell_{\bf f} (h) =
\sum_i\gamma_i  b_{\Tc_i}$, where 
$\gamma_i\in k$.	
	
		Using the expression for $\ell_f (h_j ) =
	\sum_i \gamma_{i,j} b_{j,\Tc_i}$, the above
	congruences reduce to
	$$
	\sum_{i,j} \gamma_{i,j}\beta_j b_{j,\Tc_i} \equiv 0 \mod p^N.
	$$
	The last congruence is done by
	an elementary check on the Fourier coefficients $b_{j,\Tc_i}$.

	Using {\it Krasner graded structures} in general, the abstract Kummer congruences are checked for a family of test elements (e.g. certain $p$-adic Dirichlet characters
	with values in ${\bar \Q}_p^*$, viewed as homogeneous elements of grade $j\in J\subset \Yc$, with above $J=\Yc^{alg}$.
\


	\section*{Acknowledgement}	

		Many thanks to Mirjana Vukovi\'c for the invitation to this conference, and for soliciting a paper!\\
A distinguished coauthor of Professor Marc Krasner in Paris, she went also to Moscow
in 1975-1976 to work in the Algebra Deartment MSU, with Prof. A.V.Mikhalev, E.S.Golod, Yu.I.Manin, A.I.Kostrikin I.R.Shafarevich. 	
We met there for the first that time when I was a student, then again in Moscow in 1983 and in 2009 , and in Grenoble in 2000-2001 working on graded structures and their applications, see \cite{Vu01}, talking to S.Boecherer,  Y.Varshavsky, 	R.Gillard.

\

My deep gratitude go to Khoai Ha Huy, Siegfried Boecherer and Vladimir Berkovich for fruitful discussions during a visit of the author to Vietnam in April 2015.

\





\setbeamertemplate{footline}
{%
\begin{beamercolorbox}{section in foot} 
\vskip-10pt
\fbox{
\insertpagenumber}
{}
\end{beamercolorbox}%
}
\

\ifslide
{\large\blue References}
\fi


\bibliographystyle{plain}

\begin{thebibliography}{}

\end{thebibliography}


\begin{thebibliography}{10xxxx}
\bibitem[Am75]{Am75}
{\sc Amice,  Y.}, 
{\em  Les nombres $p$-adiques}, 1975, PUF, Collection SUP, p.189.
\bibitem[Am-V]{Am-V}
{\sc Amice,  Y.} and {\sc Vélu,  J.}, 
{\em Distributions $p$-adiques associ\'ees aux
	s\'eries de Hecke}, 
Journ\'ees Arithm\'etiques de Bordeaux (Conf. Univ.
Bordeaux, 1974), Ast\'e\-risque no. 24/25, Soc. Math. France, Paris
1975, pp. 119-131



\bibitem[Boe85]{Boe85}
{\sc B\"ocherer}, S., \ {\em \"Uber die Funk\-tio\-nal\-glei\-chung 
	auto\-mor\-pher $L$--Funk\-tio\-nen zur Sie\-gel\-scher Modul\-gruppe.} 
J. reine angew. Math. 362 (1985) 146-168

\bibitem[BoeNa13]{BoeNa13}
{\sc Boecherer,} Siegfried,  {\sc  Nagaoka},  Shoyu , 
{\em On $p$-adic properties of Siegel modular forms}, arXiv:1305.0604 [math.NT]




\bibitem[Boe-Pa9]{Boe-Pa9}
{\sc B\"ocherer}, S.,  {\sc Panchishkin},  A.A.,
{\it $p$-adic Interpolation of Triple $L$-functions: Analytic Aspects}. In:
Automorphic Forms and $L$-functions II: Local Aspects -- David Ginzburg, 
Erez Lapid, 
and David Soudry, 
Editors, AMS, 2009, 313 pp.; pp.1-41

\bibitem[Boe-Pa11]{Boe-Pa11}
{\sc B\"ocherer}, S.,  {\sc Panchishkin},  A.A.,
{\it Higher Twists and Higher Gauss Sums}
Vietnam Journal of Mathematics 39:3 (2011) 309-326


\bibitem[BS00]{BS00}
{\sc B\"ocherer, S.},  and {\sc Schmidt, C.-G.}, 
{\em $p$-adic measures attached to Siegel modular forms}, 
Ann. Inst. Fourier 50, \Numero 5, 1375-1443 (2000).







\bibitem[CourPa]{CourPa}
{\sc Courtieu,}M.,  {\sc Panchishkin },A.A.,
{\em Non-Archimedean $L$-Functions and Arithmetical Siegel Modular Forms},
Lecture Notes in Mathematics 1471, Springer-Verlag, 2004 (2nd augmented ed.)



\bibitem[De84]{De84}
{\sc Deligne, P.},  {\em  Les corps locaux de caractéristique p, limites de corps locaux de caractéristique 0}. Representations of reductive groups over a local field, 119-157, Travaux en Cours, Hermann, Paris, 1984.




\bibitem[EE]{EE}
{\sc Eischen}, Ellen E.., {\em $p$-adic Differential Operators on
	Automorphic Forms on Unitary Groups.} Annales de l'Institut
Fourier 62, No.1 (2012) 177-243.

\bibitem[EHLS]{EHLS}
{\sc Eischen} Ellen E., {\sc  Harris}, Michael,  {\sc Li}, Jian-Shu,
{\sc Skinner}, Christopher M., {\em $p$-adic $L$-functions for Unitary Shimura Varieties,
	p-adic L-functions for unitary groups, part II: zeta-integral calculations},
{\tt	arXiv:1602.01776 [math.NT]}







\bibitem[GeSha]{GeSha}
{\sc Gelbart,} S.,  and  {\sc Shahidi}, F.,
{\em Analytic Properties of Automorphic $L$-functions}, Academic Press, New York, 1988.

\bibitem[GPSR]{GPSR}  {\sc Gelbart S.,Piatetski-Shapiro I.I., Rallis S.}
{\it Explicit
	constructions of automorphic $L$-functions.}
Springer-Verlag, Lect. Notes in Math.
N 1254 (1987) 152p.



\bibitem[Ha81]{Ha81}
{\sc Harris, } M., {\em
	Special values of zeta functions attached to Siegel modular
	forms. Ann. Sci.} Ecole Norm Sup. 14 (1981), 77?120. 




\bibitem[Ich12]{Ich12}
{\sc Ichikawa, T.},
{\em Vector-valued $p$-adic Siegel modular forms}, J.
reine angew. Math., DOI 10.1515/ crelle-2012-0066.

\bibitem[Ich13]{Ich13}
{\sc Ichikawa, T.},
{\em  Arithmeticity of vector-valued Siegel modular
	forms in analytic and p-adic cases}. Arxiv: 1508.03138v2  [MathNT].

\bibitem[Ka76]{Ka76}  
{\sc  Katz, N.M.}, {\em $p$-adic interpolation of real analytic Eisenstein
	series.} Ann. of Math. 104 (1976) 459--571




\bibitem[KiNa16]{KiNa16}  
{\sc Kikuta}, Toshiyuki,  {\sc Nagaoka}, Shoyu,  {\em 
	Note on mod $p$ property of Hermitian modular forms} arXiv:1601.03506 [math.NT]


\bibitem[Ko80]{Ko80}
{\sc Koblitz}, Neal, 
{\em p-adic Analysis. A Short Course on
	Recent Work}, Cambridge Univ. Press, 1980
\bibitem[Kr44]{Kr44}
{\sc Krasner}, M., {\em Une généralisation de la notion de corps-corpoïde. Un corpoïde remarquable de la
théorie des corps valués}, C. R. Acad. Sci. Paris Sér. IMath. 219 (1944), 345-347.

\bibitem[Kr50]{Kr50}
{\sc Krasner}, M., {\em Quelques méthodes nouvelles dans la théorie des corps valués complets.}
Algèbre et Théorie des Nombres, 29-39. Colloques Internationaux du Centre National
	de la Recherche Scientifique, no. 24, Centre National de la Recherche Scientifique,
	Paris, 1950.

\bibitem[Kr66]{Kr66}
{\sc Krasner}, M., Prolongement analytique uniforme et multiforme,
Collogue C.N.R.,S. n° 143, Clermont-Ferrand, 1963, Paris, Ed. C.N.R.S.,
1966, p. 97-141


\bibitem[Kr74]{Kr74}
{\sc Krasner}, M., 
{\em Rapport sur le prologement analytique dans les corps
	values complete par la methode des e1ements analytiques quasiconnexes},
Bull. Soc. Math. France, Mem. 39-40 (1974) 131-254.


\bibitem[Kr80]{Kr80}
{\sc Krasner}, M., 
{\em  Anneaux gradués généraux}, Colloque d'Algèbre Rennes (1980), 209-308.

\bibitem[KrKa51]{KrKa51}
{\sc Krasner}, M., {\sc Kaloujnine}, L.  
{\em  Produit complet des groupes de permutations et problème d'extension de groupes II }, Acta Sci. Math. Szeged , 14 (1951) p. 39-66 et 69-82.

\bibitem[KrVu87]{KrVu87}
{\sc Krasner}, M., {\sc Vukovi\'c}, M. 
{\em Structures paragraduées (groupes, anneaux, modules)}, Queen's Papers in Pure and Applied Mathematics, 77, Queen's University, Kingston, Ontario, Canada, 1987.




\bibitem[LangMF]{LangMF}
{\sc
	Lang}, Serge. {\em Introduction to modular forms. With appendixes by D. Zagier and Walter Feit.} 
Springer-Verlag, Berlin, 1995

\bibitem[Ma73]{Ma73} {\sc  Manin, Yu. I.},
{\em  Periods of cusp forms and p-adic Hecke series}, Mat.
Sbornik, 92 , 1973, pp. 378-401

\bibitem[Ma76]{Ma76} {\sc  Manin, Yu. I.},
{\em  Non-Archimedean integration and Jacquet-Langlands p-adic L-functions}, Uspekhi Mat. Nauk, 1976,
Volume 31, Issue 1(187), 5-54




\bibitem[MaPa]{MaPa}{\sc  Manin, Yu. I.},{\sc Panchishkin, A.A.},
{\em Introduction to Modern Number Theory: Fundamental Problems, Ideas and Theories} (Encyclopaedia of Mathematical Sciences), Second Edition, 504 p., Springer (2005)  



\bibitem[Pa91]{Pa91}{\sc Panchishkin, A.A.},
{\em Non-Archimedean $L$-Functions of Siegel and Hilbert Modular Forms.}
 Volume 1471 (1991)

\bibitem[PaSE]{PaSE}  {\sc Panchishkin, A.A.},
{\it On the Siegel-Eisenstein measure and
	its applications},
Israel Journal of
Mathemetics, 120, Part B (2000) 467-509.

\bibitem[PaMMJ]{PaMMJ} {\sc Panchishkin,   A.A.},
{\em A new method of constructing $p$-adic $L$-functions associated 
	with modular forms}, 
Moscow Mathematical Journal, 2 (2002), Number 
2, 1-16


\bibitem[PaTV]{PaTV} {\sc Panchishkin,   A.A.},
{\em Two variable $p$-adic $L$ functions attached to eigenfamilies of positive slope},
Invent. Math. v. 154, N3 (2003), pp. 551 - 615






\bibitem[PaIsr11]{PaIsr11}  {\sc Panchishkin, A.A.},
{\it 
	Families of Siegel modular forms, $L$-functions
	and modularity lifting conjectures.}
Israel Journal of
Mathemetics, 185 (2011), 343-368

\bibitem[Sch12]{Sch12}
{\sc Scholze,  Peter}, 
{\em Perfectoid spaces}, Publications mathématiques de l'IHÉS, vol. 116, no 1, 2012, p. 245-313

\bibitem[Se73]{Se73}
{\sc Serre,  J.--P.}, 
{\em Cours d'arithmétique}. Paris, 1970.

\bibitem[Se73]{Se73}
{\sc Serre,  J.--P.}, 
{\em Formes modulaires et fonctions z\^eta $p$-adiques}, 
Lect Notes in Math. 350 (1973) 191--268 (Springer Verlag)






\bibitem[Shi00]{Shi00}
{\sc Shimura} G.,
{\em Arithmeticity in the theory of automorphic forms},
Mathematical Surveys and Monographs, vol. 82 (Amer. Math.
Soc., Providence, 2000).


\bibitem[MC]{MC}
{\sc Skinner},  C.  and {\sc Urban}, E.  {\em
	The Iwasawa Main Cconjecture for  GL(2)}.
{\tt http://www.math.jussieu.fr/\~{}urban/eurp/MC.pdf}


\bibitem[Ta71]{Ta71}
{\sc  Tate, John},
{\em Rigid analytic spaces}, Inventiones Mathematicae, 12,  257-289 (1971)

\bibitem[U14]{U14}
{\sc  Urban,  Eric},
{\em Nearly Overconvergent Modular Forms},
Iwasawa Theory 2012.
Contributions in Mathematical and Computational Sciences Volume 7, 2014, pp 401-441, 
Date: 12 Nov 2014
{\tt http://link.springer.com/chapter/10.1007/978-3-642-55245-8$\_$14}






\bibitem[Vu01]{Vu01}
{\sc Vukovi\'c}, M.  Structures graduées et paragraduées, Prépublication de l'Institut Fourier no 536 (2001)\\
{\text https://www-fourier.ujf-grenoble.fr/sites/default/files/ref$\_$536.pdf}
\bibitem[Wa82]{Wa82}
{\sc Washington,  L.},
{\em Introduction to cyclotomic fields},
Springer
Verlag: N.Y. e.a., 1982
\bibitem[WLeo]{WLeo}
{\sc Wikipedia},
{\em Leopoldt's conjecture}.
	From Wikipedia, the free encyclopedia.



\bibitem[YS]{YS}
{\sc Yoshida, H.},   {\em Review on Goro Shimura, Arithmeticity in the
	theory of automorphic forms \cite{Shi00}}, 
Bulletin (New Series) of
the AMS, vol. 39, N3 (2002), 441-448.

\bibitem[Z13]{Z13}
{\sc Zemel,S.},   {\em 
On quasimodular forms, almost holomorphic modular forms,
and the vector-valued modular forms of Shimura.} arXiv:1307.1997
(2013)

\end{thebibliography}

\end{document}